\documentclass[12pt]{article}
\input epsf
\usepackage{amssymb}
\usepackage{latexsym}

\def \be{\begin{equation}}
\def \ee{\end{equation}}
\def \berr{\begin{eqnarray}}
\def \err{\end{eqnarray}}

\def \Om{\Omega}

\def \F{{\cal F}}

\def \U{{\cal U}_q}

\def \({\left(}
\def \){\right)}
\def \<{\langle}
\def \>{\rangle}
\def \[{\left[}
\def \]{\right]}

\def\reps{representations }
\def\rep{representation }

\newcommand{\sect}[1]{\setcounter{equation}{0}\section{#1}}

\evensidemargin \oddsidemargin
\begin{document}


\begin{center}  

{\Large\bf
Quantum Anti--de Sitter Space at Roots of Unity
 \\[4ex]}
Harold \ Steinacker 
\footnote{Harold.Steinacker@physik.uni--muenchen.de}
\\[2ex] 
{\small\it 
        Sektion Physik der Ludwig--Maximilians--Universit\"at\\
        Theresienstr.\ 37, D-80333 M\"unchen  \\[1ex] }
\end{center}

\vspace{5ex}

\begin{abstract}
An algebra of functions on q--deformed Anti-de Sitter space $AdS_q^D$ 
with star-structure is defined for roots of unity, which is
covariant under $U_q(so(2,D-1))$. The scalar fields have an intrinsic 
high--energy cutoff, and arise most naturally on 
products of the quantum AdS space with a classical sphere.
Hilbert spaces of scalar fields are constructed.
\end{abstract}

\sect{Introduction}

The aim of these notes is to give a short summary of the paper
\cite{ads}. We define and study 
a non--commutative version of the $D$--dimensional Anti--de Sitter
(AdS) space which is covariant under the standard Drinfeld--Jimbo quantum group
$U_q(so(2,D-1))$, for $q$ a root of unity and $D \geq 2$. 
The symmetry group $SO(2,D-1)$ plays the role of the
$D$--dimensional Poincar\'e group.
Part of the motivation for doing this is an interesting conjecture 
relating string or M theory 
on $AdS^{D} \times W$ with (super)conformal field theories on the boundary 
\cite{maldacena},
where $W$ is a certain sphere or a product space containing a sphere. 
There is in fact some evidence
that a full quantum treatment would lead to some non--classical 
version of the manifolds. This includes the appearance of a 
``stringy exclusion principle'' \cite{strominger} in the spectrum
of fields on AdS space.
In view of this and the well--known connections of conformal field theory 
with quantum groups at roots of unity, 
it seems plausible that some kind of quantum AdS  space at roots of unity
should be relevant to string theory in the above context.
We argue in \cite{ads} that the space presented below follows quite 
uniquely form general covariance assumptions. 

It is well--known that at roots of unity, quantum groups show 
completely new, ``non--perturbative'' features. In particular,
there exist finite--dimensional unitary \reps of 
the quantum AdS groups at roots of unity \cite{nc_reps,dobrev}, 
which consistently combine all the features
of the classical one--particle \reps with a high--energy cutoff. 
In particular, scalar fields which are unitary \reps of $U_q(so(2,D-1))$
will be obtained as polynomials in the coordinate functions.
Moreover, it turns out that the quantum spaces are obtained most naturally as
products
$AdS_q^{D} \times S^{D}/\Gamma$ if $D$ is odd, and 
$AdS_q^{D} \times S_{\chi}^{2D-1}/\Gamma$ if $D$ is even.

\sect{The structure of $AdS_q^{D}$}

The algebra of functions on quantum AdS space is a real sector of 
the complex orthogonal quantum sphere \cite{FRT}
\be
(P^-)^{ij}_{kl} t_i t_j =0, \qquad
g^{ij} t_i t_j = R^2
\label{t_2}
\ee
where $P^-$ is the $q$--antisymmetrizer \cite{FRT}.
This is covariant under $\U^{res}:= U_q^{res}(so(2,D-1))$, which has 
generators $\{H_i, X_i^{\pm}\}$ which satisfy the standard Drinfeld--Jimbo
commutation relations, plus additional generators of a {\em classical} 
Lie algebra which is either $so(D+1)$ if $D$ is odd, or $sp(D)$ if 
$D$ is even; compare \cite{nc_reps}. 
The reality condition to specify the appropriate noncompact form is 
$H_i^* = H_i$, 
${X_i^{\pm}}^*  = (-1)^E X_i^{\mp} (-1)^E = s_i X_i^{\mp}$
with $s_i = \pm 1$, where
$E$ is the generator of $U_q(so(2,D-1))$ corresponding to the energy. 
We consider roots of unity $q=e^{i \pi /M}$, and define
$M_S = M/d_S$ where $d_S = 2$ if $D$ is even, and $d_S = 1$ if $D$ is odd.

Consider the space of irreducible polynomials of degree $n$ in the 
generators, $\F(n) =  \U^{res} \cdot (t_1)^n$. This is an irreducible 
highest--weight \rep of $\U^{res}$, and it is unitary with respect to
the compact form $U_q^{res}(so(D+1))$ if 
\be
n < k_S:= M_S - (D-2)/2.
\ee
Hence in that case, $\F(n)$ should be considered as a space of functions
on the compact quantum sphere $S_q^{D-1}$. 
However if $M_S \leq n < M_S + k_S$, then $\U^{fin} \cdot (t_1)^n$
is an irreducible unitary \rep of $\U^{fin} = U_q^{fin}(so(2,D-1))$, 
and can be identified
with a square--integrable scalar field with the same lowest weight (i.e. mass)
on the classical AdS space. Here $\U^{fin}$ is the subalgebra of 
$\U^{res}$ with generators $\{H_i, X_i^{\pm}\}$, 
which becomes the usual AdS group in the classical limit.
Therefore it make sense to define the Hilbert space of 
(positive--energy, square--integrable) functions on quantum AdS space by
\be
AdS_q^{D}:= \bigoplus_{M_S \leq n < M_S + k_S}\;\; \U^{fin} \cdot (t_1)^n.
\label{real_AdS}
\ee
From a physical point of view, this describes spin 0 particles. 

The coordinate functions $t_i$ define operators $\hat t_i$ 
on this 
Hilbert space, which induces a star structure on them by the
operator adjoint. It can be written explicitly in the form 
$\hat t_i^* = - \Om (-1)^E \; \hat t_j C^j_i \Om (-1)^E$, where
$\Om$ is a unitary element with $\Om^2=1$, which implements the 
the longest element of the Weyl group.

Upon closer examination, it turns out that the space of (``almost'') all
polynomial functions in the generators $t_i$ decomposes as 
\be
\oplus_n \F(n) = 
\Big( AdS_q^{D}\times M/\Gamma \Big)  \;\; \oplus \;\;  
\Big( S_q^{D} \times M/\Gamma \Big)  \;\; \oplus \dots,
\ee
where $M = S^{D}$ or $S_{\chi}^{2D-1}$, and $M/\Gamma$ 
denotes the space of functions on some (twisted)
orbifold; see \cite{ads} for more details.
This means that the space which arises most naturally is a 
product of quantum AdS space 
with an orbifold of a classical sphere or of a ``chiral'' sphere, plus
other sectors. The symmetries of these additional spheres are implemented 
by the classical generators which are contained in the full 
algebra $U_q^{res}(so(2,D-1))$, as mentioned above.

One can also look at this additional structure from an other point of view. 
For example in the case $D=4$, it turns out that the smallest irreducible 
\rep of $U_q^{res}(so(2,3))$ which contains a scalar field as in 
(\ref{real_AdS}) on quantum 
$AdS$ space is in fact a 4--dimensional multiplet of
$sp(4)$, which is in a sense spontaneously broken to
$su(2) \times su(2)$. In other words, there are in general additional
``global'' symmetries, and the sum of all these multiplets in $\oplus \F(n)$
becomes a sector of a classical sphere.
 
Furthermore, we give an argument in \cite{ads} that this quantum AdS space
has an intrinsic length scale of order $L_{min} = R/M_S$, where
the geometry is expected to become non--classical. 
One can also define a differential calculus and integration.
A formulation of field theory on this space
should hence be possible along the lines of \cite{fuzzy_q}.

\end{document}